\documentclass[12pt]{article}
\usepackage{oldlfont}
\usepackage{enumerate}
\usepackage{amsmath}
\usepackage{amssymb}
\usepackage{amsthm}
\usepackage{mathrsfs}
\usepackage{amscd}
\usepackage{exscale}
\usepackage{latexsym}
\usepackage[all]{xy}
\usepackage{hyperref}
\usepackage{layout}
\usepackage{color}
\usepackage{graphicx}
\DeclareFontFamily{U}{wncy}{}
\DeclareFontShape{U}{wncy}{m}{n}{<->wncyr10}{}
\DeclareSymbolFont{mcy}{U}{wncy}{m}{n}
\DeclareMathSymbol{\Sh}{\mathord}{mcy}{"58} 

\begin{document}

\baselineskip=17pt

\pagestyle{headings}

\numberwithin{equation}{section}

\makeatletter                                                           

\def\section{\@startsection {section}{1}{\z@}{-5.5ex plus -.5ex         
minus -.2ex}{1ex plus .2ex}{\large \bf}}                                 


\renewcommand{\sectionmark}[1]{\markboth{ #1}{ #1}}
\renewcommand{\subsectionmark}[1]{\markright{ #1}}

\addtolength{\headheight}{0.5pt} 

\newtheorem{thm}{Theorem}
\newtheorem{mainthm}[thm]{Main Theorem}

\newcommand{\ZZ}{{\mathbb Z}}
\newcommand{\GG}{{\mathbb G}}
\newcommand{\Z}{{\mathbb Z}}
\newcommand{\RR}{{\mathbb R}}
\newcommand{\NN}{{\mathbb N}}
\newcommand{\GF}{{\rm GF}}
\newcommand{\QQ}{{\mathbb Q}}
\newcommand{\CC}{{\mathbb C}}
\newcommand{\FF}{{\mathbb F}}

\newtheorem{lem}[thm]{Lemma}
\newtheorem{con}[thm]{Conjecture}
\newtheorem{cor}[thm]{Corollary}
\newtheorem{pro}[thm]{Proposition}
\newtheorem{proprieta}[thm]{Property}
\newcommand{\pf}{\noindent \textbf{Proof.} \ }
\newcommand{\eop}{$_{\Box}$  \relax}
\newtheorem{obs}[thm]{Remark}
\newtheorem{num}{equation}{}

\theoremstyle{definition}
\newtheorem{rem}[thm]{Remark}
\newtheorem*{D}{Definition}

\newcommand{\nsplit}{\cdot}
\newcommand{\G}{{\mathfrak g}}
\newcommand{\GL}{{\rm GL}}
\newcommand{\SL}{{\rm SL}}
\newcommand{\SP}{{\rm Sp}}
\newcommand{\LL}{{\rm L}}
\newcommand{\Ker}{{\rm Ker}}
\newcommand{\la}{\langle}
\newcommand{\ra}{\rangle}
\newcommand{\PSp}{{\rm PSp}}
\newcommand{\U}{{\rm U}}
\newcommand{\GU}{{\rm GU}}
\newcommand{\Aut}{{\rm Aut}}
\newcommand{\Alt}{{\rm Alt}}
\newcommand{\Sym}{{\rm Sym}}

\newcommand{\isom}{{\cong}}
\newcommand{\z}{{\zeta}}
\newcommand{\Gal}{{\rm Gal}}

\newcommand{\F}{{\mathbb F}}
\renewcommand{\O}{{\cal O}}
\newcommand{\Q}{{\mathbb Q}}
\newcommand{\R}{{\mathbb R}}
\newcommand{\N}{{\mathbb N}}
\newcommand{\A}{{\mathcal{A}}}
\newcommand{\B}{{\mathcal{B}}}
\newcommand{\E}{{\mathcal{E}}}
\newcommand{\J}{{\mathcal{J}}}


\newcommand{\DIM}{{\smallskip\noindent{\bf Proof.}\quad}}
\newcommand{\CVD}{\begin{flushright}$\square$\end{flushright}
\vskip 0.2cm\goodbreak}


\vskip 0.5cm

\title{On the local-global divisibility over abelian varieties}
\author{Florence Gillibert\thanks{supported by the Fondecyt Iniciaci\'on Project 11130409}, Gabriele Ranieri\thanks{supported by the Fondecyt Regular Project 1140946}}
\date{  }
\maketitle

\vskip 1.5cm

\begin{abstract}
Let $p \geq 2$ be a prime number and let $k$ be a number field. 
Let $\mathcal{A}$ be an abelian variety defined over~$k$.
We prove that if $\Gal ( k ( \A[p] ) / k )$ contains an element $g$ of order dividing $p-1$ not fixing any non-trivial element of $\A[p]$ and $H^1 ( \Gal ( k ( \A[p] ) / k ), \A[p] )$ is trivial, then the local-global divisibility by $p^n$ holds for $\A ( k )$ for every $n \in \N$. Moreover, we prove a similar result without the hypothesis on the triviality of $H^1 ( \Gal ( k ( \A[p] ) / k ) , \A[p] )$, in the particular case where $\A$ is a principally polarized abelian variety.
Then, we get a more precise result in the case when $\A$ has dimension~$2$.
Finally we show with a counterexample that the hypothesis over the order of~$g$ is necessary. 

In the Appendix, we explain how our results are related to a question of Cassels on the divisibility of the Tate-Shafarevich group, studied by Ciperani and Stix \cite{C-S}. 

\end{abstract}

\section{Introduction}

Let $k$ be a number field and let ${\mathcal{A}}$ be a commutative algebraic group defined over $ k $.
Several papers have been written on  the following classical question, known as \emph{Local-Global Divisibility Problem}.

\par\bigskip\noindent  P{\small ROBLEM}: \emph{Let $P \in {\mathcal{A}}( k )$. Assume that for all but finitely many valuations $v$ of $k$, there exists $D_v \in {\mathcal{A}}( k_v )$ such that $P = qD_v$, where $q$ is a positive integer. Is it possible to conclude that there exists $D\in {\mathcal{A}}( k )$ such that $P=qD$?}

\par\bigskip\noindent  By  B\'{e}zout's identity, to get answers for a general integer it is sufficient to solve it for powers $p^n$ of a prime. In the classical case of ${\mathcal{A}}={\mathbb{G}}_m$ and $k = \Q$, the answer is positive for $ p $ odd, and negative for instance for $q=8$ (and $P=16$) (see for example \cite{AT}, \cite{T}).

\bigskip  For general commutative algebraic groups, Dvornicich and Zannier gave a cohomological interpretation of the problem (see \cite{DZ} and \cite{DZ3}) that we shall explain.
Let $ \Gamma $ be a group and let $ M $ be a $\Gamma$-module.
We say that a cocycle $Z \colon \Gamma \rightarrow M$ satisfies the local conditions if for every $\gamma \in \Gamma$ there exists $m_\gamma \in M$ such that $Z_\gamma = \gamma ( m_\gamma ) - m_\gamma$.
The set of the class of cocycles in $ H^1 ( \Gamma , M ) $ that satisfy the local conditions is a subgroup of $H^1 ( \Gamma, M )$.
We call it the first local cohomology group $H^1_{{\rm loc}} ( \Gamma, M )$.
Equivalently,
\[
H^1_{{\rm loc}} ( \Gamma , M ) = \cap_{ C \leq \Gamma } \ker ( H^1 ( \Gamma , M ) \rightarrow H^1 ( C , M ) ),
\]
where $C$ varies among the cyclic subgroups of $ \Gamma $ and the above maps are the restrictions.   
Dvornicich and Zannier \cite[Proposion 2.1]{DZ} proved the following result.

\begin{pro}\label{pro1}
Let $ p $ be a prime number, let $n$ be a positive integer, let $ k $ be a number field and let $ \A $ be a commutative algebraic group defined over $ k $.
If $H^1_{{\rm loc}} ( \Gal ( k ( \A[p^n] ) / k ) , \A[p^n] ) = 0$, then the local-global divisibility by $p^n$ over $ \A ( k ) $ holds.
\end{pro}

The converse of Proposition \ref{pro1} is not true.
However, in the case when the group $H^1_{{\rm loc}} ( \Gal ( k ( \A[p^n] ) / k ) , \A[p^n] )$ is not trivial we can find an extension $ L $ of $k$ $k$-linearly disjoint with $k ( \A[p^n] )$ in which the local-global divisibility by $ p^n $ over $ \A ( L ) $ does not hold (see \cite[Theorem 3]{DZ3} for the details).

From now on let $ p $ be a prime number, let $ k $ be a number field and let $ \E $  be an elliptic curve defined over $ k $.   
Dvornicich and Zannier (see \cite[Theorem 1]{DZ3}) found a geometric criterion for the validity of the local-global divisibility principle by a power of $ p $ over $ \E ( k ) $.
In \cite{P-R-V} and \cite{P-R-V2} Paladino, Viada and the second author refined this criterion.
Ciperani and Stix  \cite[Theorem A, Theorem B]{C-S} also proved a similar criterion to give an answer to a question of Cassels on elliptic curves (see the Appendix).  
Moreover, very recently, T. Lawson and C. Wuthrich \cite{LW} obtained a very strong criterion for the vanishing of the first cohomology group of the Galois module of the torsion points of an elliptic curve defined over $ \Q $ that allowed them to find a simplest proof of the main result of \cite{P-R-V2}.
They also found a counterexample to the local-global divisibility by powers of $3$ for an elliptic curve over $ \Q $. 
From this result, the examples of Dvornicich and Zannier \cite{DZ2} and Paladino \cite{Pal}, \cite{Pal2} and the main result of \cite{P-R-V2}, it follows that the set of prime numbers $q$ for which there exists an elliptic curve $\E^\prime$ defined over $ \Q $ and $n \in \N$ such that the local-global divisibility by $q^n$ does not hold over $\E^\prime ( \Q )$ is just $\{2, 3\}$.
 
Let us consider now an arbitrary abelian variety. To our knwowledge the unique known geometric criterion for the validity of the local-global divisibility principle by a power of $ p $ for an abelian variety of dimension $>1$ over a number field was proved by Ciperani and Stix (see \cite[Theorem D]{C-S}. For a connection with this result and the local-global divisibility problem see \cite[Remark 20]{C-S} and the Appendix).

The results on elliptic curves and this last result gave a motivation to look for other geometric criterions for the local-global divisibility principle over the family of the abelian varieties. 
From now on, let $ \A $ be an abelian variety defined over $ k $ of dimension $d \in \N^\ast$. 
Moreover, for every positive integer $ n $, we set $K_n = k ( \A[p^n] )$ and $G_n = \Gal ( K_n / k )$.
We prove the following result.

\begin{thm}\label{teo2}
Suppose that $ G_1 $  contains an element $ g $ whose order divides $ p-1 $ and not fixing any non-trivial element of $ \A[p] $. Moreover suppose that $H^1 ( G_1, \A[p] ) = 0$. Then $H^1_{{\rm loc}} ( G_n , \A[p^n] ) = 0$ for every positive integer $n$ and so the local-global divisibility by $p^n$ holds for $\A ( k )$.
\end{thm}

Let us now fix a polarization on $ \A $ over $ k $ and let us suppose that $ p $ does not divide the degree of the polarization. 
We prove the following result, in which there is no hypothesis on $ H^1 ( G_1 , \A[p] ) $.
However, we need a hypothesis on the field $ k $.

\begin{thm}\label{teo5}   
Let $ \A $ be a polarized abelian variety of dimension $ d $ defined over $ k $ and let $ p $ be a prime not dividing the degree of the polarization. 
Suppose that $k \cap \Q ( \zeta_p ) = \Q$.
Set $i = ( ( 2d )! , p-1 )$ and $ k_i $ the subfield of $k ( \zeta_p )$ of degree $i$ over $ k $.
If for every non-zero $P \in \A[p]$ the field $k ( P ) \cap k ( \zeta_p )$ strictly contains $k_i$, then for every positive integer $n$, the group $H^1_{{\rm loc}} ( G_n , \A[p^n] ) = 0$ and so the local-global divisibility by $p^n$ holds for $\A ( k )$.
\end{thm}  

Suppose now that $ \A $ has dimension $2$.
By using Theorems \ref{teo2}, \ref{teo5} and the results on sections \ref{sec2} and \ref{sec12}, we shall give a much more precise criterion, which is a weak generalization to abelian surfaces of the main result of \cite{P-R-V} on elliptic curves.

\begin{thm}\label{teo6}   
Let $ \A $ be a polarized abelian surface defined over $ k $.
For every prime number $p > 3840$ such that $k \cap \Q ( \zeta_p ) = k$ and not dividing the degree of the polarization, if there exists $n \in \N$ such that $H^1_{{\rm loc}} ( G_n , \A[p^n] ) \neq 0$, then there exists a finite extension $ \widetilde{k} $ of $ k $ of degree $\leq 24$ over $k$ such that $ \A $ is $\widetilde{k}$-isogenous to an abelian surface with a torsion point of order $ p $ defined over $ \widetilde{k} $.  
\end{thm} 

Merel \cite{M} made the following conjecture on the torsion of abelian varieties over a number field and proved it in the case of dimension $1$.

\begin{con}\label{con1}[Merel's Conjecture]
Let $ d $ and $m$ be positive integers. There exists a positive constant $C ( d, m )$, only depending on $d$ and $m$, such that for every abelian variety of dimension $ d $ defined over a number field $ k $ of degree $m$, and for every prime number $p > C ( d, m )$, $\A$ does not admit any point of order $ p $ defined over $ k $.
\end{con}

Then, we have the following corollary of Theorem \ref{teo6}.

\begin{cor}\label{cor7Merel}
If Merel's Conjecture is true, then for every positive integer $ m $ there exists a constant $C ( m )$ only depending on $m$ such that for every principally polarized abelian surface $ \A $ defined over a number field $ k $ of degree $m$ over $\Q$ and for every prime number $p > C ( m )$, then for every positive integer $n$ the local-global divisibility by $p^n$ holds for $\A ( k )$.
\end{cor}

Here it is the plan of this paper.
In Section \ref{sec2} we prove some algebraic results necessary for the proof of Theorem \ref{teo2}, Theorem \ref{teo5} and Theorem \ref{teo6}.
Moreover we prove Theorem \ref{teo2}.

For every prime number $ p $ not dividing the degree of a polarization of a polarized abelian variety, the image of the absolute Galois group on the group of the automorphism of the $p$-torsion is contained in the group of the symplectic similitudes for the Weil-pairing.  
In Section \ref{sec12} we describe such a group and we prove Theorem \ref{teo5}.  
For the proof of Theorem \ref{teo6} is necessary a very precise study of the properties of the group $ {\rm GSp}_4 ( \F_p ) $.
We do this in section \ref{sec3} and then we finish such section by proving Theorem \ref{teo6}.
In Section \ref{sec6} we give an example that shows that the hypothesis on the order of $ g $ in Theorem \ref{teo2} is necessary.
Finally we explain in an Appendix the connection with the local-global divisibility problem and a question of Cassels studied in particular by Ciperani and Stix in \cite{C-S}.          

\section{Algebraic preliminaries}\label{sec2}

\subsection{Coprime groups and cohomology}

Classical Frattini's theory (see for instance \cite{A}) is very useful to prove the following Proposition, which is the first step to prove Theorem \ref{teo2}.

\begin{pro}\label{pro11}
Let $p$ be a prime number and let $ G $ be a finite group such that $G = \langle g, H \rangle$, where $ g $ has order dividing $p-1$ and $H$ is a $p$-group, which is normal in $ G $.
There exist $r \in \mathbb{N}$ and a generator set $\{h_1, h_2, \ldots , h_r \}$ of $ H $ such that, for every $1 \leq i \leq r$, there exists $\lambda_i \in \Z$ such that
\[
g h_i g^{-1} = h_i^{\lambda_i}.
\]
\end{pro}

\DIM Suppose $\vert H \vert = p^m$ with $m \in \N$. 
The proof is by induction on $ m $.

If $m = 1$ we have that $H$ is cyclic generated by an element $h_1$.
Since $H$ is normal in $ G $, we have $g h_1 g^{-1} = h_1^{\lambda_1}$ for a $\lambda_1 \in \mathbb{Z}$ and there is nothing to prove.

Suppose that the assumption is true for every natural number $j < m$.
We recall that the Frattini subgroup (see \cite[p. 105]{A} for the details) $\phi ( H )$ of $H$ is the intersection of all maximal subgroups of $H$ and that $H / \phi ( H )$ is elementary abelian (i.e. is isomorphic to a finite product of groups isomorphic to $\Z / p \Z$).

Let us show that $\phi ( H )$ is normal in $ G $.
Let $M$ be a maximal subgroup of $H$.  
We have $g M g^{-1} \subseteq H$ because $H$ is normal.
Then the action by conjugation of $g$ permutes the maximal subgroups of $H$.
Then, since $\phi ( H )$ is the intersection of every maximal subgroup of $H$, it is normal in $ G $.   

We use the following well-known result.

\begin{thm} [Burnside basis theorem]
Let $H$ be a finite $p$-group.
A subset of $H$ is a set of generators for $H$ if and only if its image in $H / \phi ( H )$ is a set of generators for $H / \phi ( H )$.
\end{thm}

Consider $H / \phi ( H )$. 
Since $\phi ( H )$ is normal in $ G $ and $H / \phi ( H )$ is abelian, the function 
\[
f \colon H / \phi ( H ) \rightarrow H / \phi ( H )
\]
that sends $h \phi ( H )$ to $g h g^{-1} \phi ( H )$ is well-defined and it is actually a $\Z / p \Z$-linear isomorphism.
Since $g$ has order dividing $p-1$, also the order of $f$ divides $p-1$ and so $f$ is diagonalizable on the $\Z / p \Z$-vector space $H / \phi ( H )$. 
Then there exist $v_1, v_2, \ldots, v_k \in H$ such that $\{ v_i \phi ( H ) \ 1 \leq i \leq k \}$ is a $\Z / p \Z$-basis of $H/ \phi ( H )$ and there exist $\lambda_i \in \Z$ such that $g v_i g^{-1} \phi ( H ) = v_i^{\lambda_i} \phi ( H )$.

Suppose that $k = 1$.
Then $H / \phi ( H )$ has a unique generator $v_1 \phi ( H )$.
By Burnside basis theorem, $H$ is then generated by $v_1$ and it is cyclic. 
Since $H$ is normal in $G$, we have that $g v_1 g^{-1} = v_1^{\lambda}$ for a $\lambda \in \Z$, which is the thesis.

Suppose $k > 1$.
Consider the two groups $H_1$, $H_2 \subseteq H$, such that
\[
H_1 = \langle v_1, \phi ( H ) \rangle, \hspace{5 pt} H_2 = \langle v_2, v_3, \ldots, v_k, \phi ( H ) \rangle.
\]          
Then set $\Gamma_1$ the subgroup of $G$ generated by $g$ and $H_1$ and $\Gamma_2$ the subgroup of $G$ generated by $g$ and $H_2$.
We remark that $H_1$ is normal in $\Gamma_1$ and $H_2$ is normal in $\Gamma_2$.
In fact, as $H_1 / \phi ( H )$ is generated by $v_1$, all element of $H_1$ is in $v_1^a \phi ( H )$ for some integer $a$.   
In the same way we can prove that $H_2$ is normal in $\Gamma_2$. 

We now prove that $\Gamma_1$ and $\Gamma_2$ are not $G$.
Since $H_1$ and $H_2$ are respectively normal over $\Gamma_1$ and $\Gamma_2$ and $\Gamma_1$ and $\Gamma_2$ are generated by such groups and an element of order not divided by $p$, $H_1$ is the unique $p$-Sylow subgroup of $\Gamma_1$ and $H_2$ is the unique $p$-Sylow subgroup of $\Gamma_2$.
Since $H_1$ and $H_2$ are properly contained in $H$, we have that $\Gamma_1$ and $\Gamma_2$ are properly contained in $G$.

Then we can apply the inductive hypothesis to $\Gamma_1$ and $\Gamma_2$.
Since $H$ is generated by $H_1$ and $H_2$, a union of a set of generators of $H_1$ with a set of generators of $H_2$ gives a set of generators of $H$. 
This concludes the proof.
\CVD

The following corollary relates Proposition \ref{pro11} with the vanishing of the first local cohomology group.  

\begin{cor}\label{cor12}
Let $V_{n, d}$ be the group $( \Z/ p^n \Z )^{2d}$ and let $G$ be a subgroup of ${\rm GL}_{2d} ( \Z / p^n \Z )$ acting on $V_{n, d}$ in the usual way. Suppose that the normalizer of a $p$-Sylow subgroup $H$ of $G$ contains an element $g$ of order dividing $p-1$ such that $g-Id$ is bijective. Then $H^1_{{\rm loc}} ( G, V_{n, d} ) = 0$.
\end{cor}

\DIM Consider the two restrictions
\[
H^1 ( G, V_{n, d} ) \rightarrow H^1 ( \langle g, H \rangle , V_{n, d} ) \rightarrow H^1 ( H, V_{n, d} ).
\]      
Notice $H^1 ( G, V_{n, d} ) \rightarrow H^1 ( H, V_{n, d} )$ is injective since $V_{n, d}$ is a $p$-group, and $H$ a $p$-Sylow subgroup of $G$. We deduce that $H^1 ( G, V_{n, d} ) \rightarrow H^1 ( \langle g, H \rangle , V_{n, d} )$ is also injective.
Moreover such maps induce maps on the first local cohomology group.
Then the restriction $H^1_{{\rm loc}} ( G, V_{n, d} ) \rightarrow H^1_{{\rm loc}} ( \langle g, H \rangle , V_{n, d} )$ is injective and so, to prove the corollary, it is sufficient to prove that $H^1_{{\rm loc}} ( \langle g, H \rangle , V_{n, d} ) = 0$.
Apply Proposition \ref{pro11} to $\langle g, H \rangle$.
Then there exists $r \in \N$ and generators $h_1, h_2, \ldots , h_r$ of $H$ such that, for every $1 \leq i \leq r$, $g h_i g^{-1}$ is a power of $h_i$.
For every $i$ between $1$ and $r$, set $\Gamma_i = \langle g, h_i \rangle$ and $H_i$ the cyclic group generated by $h_i$. 
For every $1 \leq i \leq r$, $H_i$ is the $p$-Sylow subgroup of $\Gamma_i$.
Then we have that $H^1_{{\rm loc}} ( \Gamma_i, V_{n, d} ) \rightarrow H^1_{{\rm loc}} ( H_i , V_{n, d} )$ is injective. 
Moreover, since $H_i$ is cyclic, $H^1_{{\rm loc}} ( H_i , V_{n, d} ) = 0$, and so $H^1_{{\rm loc}} ( \Gamma_i, V_{n, d} ) = 0$. 
Then, if $Z$ is a cocycle of $\langle g, H \rangle$ satisfying the local conditions, for every $i$ between $1$ and $r$, it is a coboundary over $\Gamma_i$ and so, for every $\gamma_i \in \Gamma_i$, there exists $v_i \in V_{n, d}$ such that $Z_{\gamma_i} = \gamma_i ( v_i ) - v_i$.
Since $g \in \Gamma_i$ for every $1 \leq i \leq r$, for every $i, j$ we have $g ( v_i ) - v_i = g ( v_j ) - v_j$.
The injectivity of $g - Id$ implies that for every $i, j$ $v_i = v_j$.
Since $g, h_1, h_2, \ldots, h_r$ generate $\langle g, H \rangle$ we get that $Z$ is a coboundary over $\langle g, H \rangle$.
\CVD

\begin{obs}\label{rem11}
Let $N$ be the normal subgroup of $G$ of the elements congruent to the identity modulo $p$ (here we use the notation of Corollary \ref{cor12}).
We shall prove (see Lemma \ref{lem13}) that if there exists $\widetilde{g}$ in $G_1=G / N$ such that $\widetilde{g}$ is in the normalizer of a $p$-Sylow, with order dividing $p-1$ and such that $\widetilde{g} - Id$ is bijective, then there exists $g \in G$ in the normalizer of a $p$-Sylow with order diving $p-1$ and such that $g - Id$ is bijective.
\end{obs}

The existence of an element $h\in G$ such that $h  \widetilde{g}^{-1} \in N$ and $h$ in the normalizer of a $p$-Sylow of $G$ comes from the Lemma \ref{lem13} below. 
By elevating $h$ to an adequate power of $p$ we find an element $g$ compling the conditions. 

\begin{lem}\label{lem13}
Let $G$ be a group, let $N$ be a normal subgroup of $G$ and let $H$ be a $p$-Sylow of $G$. Let $g$ be an element of $g$ such that its class in $G/N$ is in the normalizer of the p-Sylow $HN/N$ of $G/N$. Then there exists an element of the class $gN$, which is in the normalizer of $H$.

In particular, if the $p$-Sylow $H$ is contained in a normal subgroup $N$ of $G$, then for every class of $G/N$, there exists an element of the class which is in the normalizer of $H$. 
\end{lem}

\DIM Let $\widetilde{g}$ be the class of $g$ modulo $N$. By hypothesis $\widetilde{g} ( HN/N ) \widetilde{g}^{-1}= HN/N$. 
We deduce that $gHNg^-1/N=HN/N$, then $gHg^{-1}N=HN$. 
So $gHg^{-1}$ and $H$ are two $p$-Sylow subgroups of $HN$.
Then there are conjugated by some element $x$ of $HN$. There exists $h\in H$ and $n\in N$ such that $x=nh$. So $gHg^{-1}=nhH(nh)^{-1}$ from which we deduce that $n^{-1}g$ is in the class $gN$ and in the normalizer of $H$.
\CVD

Lemma \ref{lem13} does not give a precise information on the order of the elements of the normalizer of $H$. 
Neverthless if $H$ is contained in a normal subgroup $N$ and $N$ and $G/ N$ have coprime orders, we have a coprime action (see \cite[Chapter 8]{A}) and so there exists a subgroup of $G$ isomorphic to $G/ N$ and disjoint to $N$.
Then in this case the normalizer contains a group isomorphic to $G/ N$.
Next corollary treats the case when $( \vert G/N \vert, \vert N \vert )$ is small (a sort of near coprime action) and it is crucial for proving Theorem \ref{teo5}.

\begin{cor}\label{cor14}    
Let $V_d$ be the group $( \Z / p \Z )^{2d}$ and let $G$ be a subgroup of the group ${\rm GL} ( 2d, \Z / p \Z )$ acting on $V_d$ in the usual way. Let $H$ be a $p$-Sylow subgroup of $G$. Suppose that there exists $N$ a normal subgroup of $G$ such that $G/ N$ is isomorphic to $\Z / ( p-1 ) \Z$. Let $i = ( ( 2d )!, p-1 )$. Then the normalizer of $H$ contains an element of order $p-1$ whose class modulo $N$ has order divided by $( p-1 ) / i$.
\end{cor}

\DIM Since $\vert G/ N \vert$ is not divisible by $p$, it is clear that $H \subseteq N$.
Let $g$ be in $G$ such that the class of $g$ modulo $N$ is a generator of $G/ N$.
By Lemma \ref{lem13} there exists an element in the class of $g$ modulo $N$ (we call it $g$ by abuse of notation) such that $g$ is in the normalizer of $H$. 
Since the class of $g$ has order $p-1$ in $G/ N$, the order of $g$ is $( p-1 ) r$, where $r$ is a positive integer.
Then $g^r$ is in the normalizer of $H$, has order $p-1$ and the order of its class in $G/ N$ is $( p-1 )/ ( p-1, r )$.
Then it is sufficient to prove that $( p-1 , r )$ divides $i$.
Since $p$ does not divide $p-1$ we can suppose that $r$ is not divisible by $p$.
Then the action of $g$ is semisimple and so there exists $c \in \N$ such that a matrix associated to $g$ can be decomposed in $c$ blocks of matrices $l_j \times l_j$ acting irreducibly over a sub-space of $V_d$, such that $\sum_{j = 1} l_j = 2d$ and the order of the $j$th block dividing $p^{l_j}-1$. 
Then the order of $g$ is the least common multiple of the order of the blocks and so $r$ divides the least common multiple of the $( p^{l_j} -1 ) / ( p - 1 )$.
Observe that 
\[
\frac {p^{l_j} - 1} { p - 1} = p^{l_j-1} + p^{l_j-2} + \ldots + 1. 
\]
We have
\[   
p^{l_j-1} + p^{l_j-2} + \ldots + 1 - ( p - 1 ) ( p^{l_j-2} + 2 p^{l_j-3} + \ldots +  ( l_j-1 ) ) = l_j.
\]   
Then $( p^{l_j} - 1 ) / ( p-1 ), ( p - 1) ) = ( l_j , p-1 )$.
Since for every $j$ $l_j \leq 2d$, the least common multiple of the $( l_j, ( p - 1) )$ divides $( 2d )!$, which proves the lemma.
\CVD

\subsection{Cocycles satisfying the local conditions and cohomology of the $p$-torsion}

For every $r$ between $1$ and $n$, let $V_{r, d}$ be the group $( \Z / p^r \Z )^d$. 

\begin{lem}\label{obs51}
Let $G$ be a subgroup of ${\rm GL}_{2d} ( \Z / p^n \Z )$ acting on $V_{n, d}$ in the classical way. Suppose that $G$ contains an element $\delta$ such that $\delta -Id$ is a bijective automorphism of $V_{n, d}$. Then the homomorphism $H^1 ( G , V_{n, d}[p] ) \rightarrow H^1 ( G , V_{n, d} )$ induced by the exact sequence of $G$-modules 
\[
0 \rightarrow V_{n, d}[p] \rightarrow V_{n, d} \stackrel{p} \rightarrow V_{n, d}[p^{n-1}] \rightarrow 0,
\]
is injective and its image is $H^1 ( G , V_{n, d} )[p]$.
In other words it induces an isomorphism between $H^1 ( G , V_{n, d}[p] )$ and $H^1 ( G , V_{n, d} )[p]$.   
\end{lem}

\DIM The following exact sequence of $G$-modules 
\[
0 \rightarrow V_{n, d}[p] \rightarrow V_{n, d} \stackrel{p} \rightarrow V_{n, d}[p^{n-1}] \rightarrow 0
\]
(here the first map is the inclusion and the second map is the multiplication by $p$) induces a long exact sequence of cohomology groups:
\[
\ldots \rightarrow H^0 ( G , V_{n, d}[p^{n-1}] ) \rightarrow H^1 ( G , V_{n, d}[p] ) \rightarrow H^1 ( G , V_{n, d} ) \rightarrow H^1 ( G, V_{n, d}[p^{n-1}] ).
\]
Since $G$ contains an element $\delta$ such that $\delta-Id$ is bijective over $V_{n, d}$, then $H^0 ( G , V_{n, d}[p^{n-1}] ) = 0$.
Hence we have the exact sequence 
\begin{equation}\label{eqn:relnuova1} 
0 \rightarrow H^1 ( G , V_{n, d}[p] ) \rightarrow H^1 ( G , V_{n, d} ) \rightarrow H^1 ( G, V_{n, d}[p^{n-1}] ).
\end{equation}
In particular  $H^1 ( G , V_{n, d}[p] ) \rightarrow H^1 ( G , V_{n, d} )$ is injective.
 
Let $Z$ be a cocycle from $G$ to $V_{n, d}$, representing a class $[Z]$ in $H^1 ( G , V_{n, d} )$ of order $p$. 
Then there exists $v \in V_{n, d}$ such that $p Z_\sigma = \sigma ( v ) - v$ for every $\sigma \in G$.
Since there exists $\delta \in G$ such that $\delta - Id$ is bijective over $V_{n, d}$ and $p Z_\delta = \delta ( v ) - v \in V_{n, d}[p^{n-1}]$, we get that $v \in V_{n, d}[p^{n-1}]$.
Then the cocycle from $G$ to $V_{n, d}[p^{n-1}]$ sending $\sigma$ to $p Z_\sigma$ for every $\sigma \in G$ is a coboundary and so the image of $[Z]$ over $H^1 ( G, V_{n, d}[p^{n-1}] )$ is $0$.  
Then there exists $[W] \in H^1 ( G , V_{n, d}[p] )$ such that the image of $[W]$ by $H^1 ( G , V_{n, d}[p] ) \rightarrow H^1 ( G , V_{n, d} )$ is $[Z]$ (see the sequence (\ref{eqn:relnuova1})).
This proves that $H^1 ( G , V_{n, d}[p] ) \rightarrow H^1 ( G , V_{n, d} )[p]$ is surjective and it concludes the proof.
\CVD

Next lemma gives the key step to prove Therorem \ref{teo2} and it will be very useful to study the local-global divisibility problem on abelian surfaces.

\begin{lem}\label{lem21}
Let $G$ be a subgroup of ${\rm GL}_{2d} ( \Z / p^n \Z )$ acting on $V_{n, d}$ in the classical way and let $H$ be the normal subgroup of $G$ of the elements acting like the identity over $V_{n, d}[p]$. Suppose that $G$ contains an element $\delta$ such that $\delta -Id$ is a bijective automorphism of $V_{n, d}$. Let $Z \colon G \rightarrow V_{n, d}$ be a cocycle whose restriction to $H$ is a coboundary. If $H^1 ( G/ H , V_{n, d}[p] ) = 0$, then $Z$ is a coboundary.
\end{lem}

\DIM Consider the following commutative diagram:
$$
\begin{array}{ccccccc}
0 & \rightarrow & H^1 ( G , V_{n, d}[p] ) & \rightarrow & H^1 ( G , V_{n, d} )[p] &\rightarrow & 0 \\
 &  & \downarrow  & & \downarrow & &   \\
 0 & \rightarrow & H^1 ( \langle \delta , H  \rangle , V_{n, d}[p] ) & \rightarrow & H^1 ( \langle \delta , H \rangle , V_{n, d} )[p]  & \rightarrow & 0,
\end{array}
$$
where the isomorphisms on the lines are the functions of Lemma \ref{obs51}, and the functions on the columns are the restrictions.
Then, since the restriction $H^1 ( \langle \delta , H  \rangle , V_{n, d} ) \rightarrow H^1 ( H , V_{n, d} )$ is injective, if $\ker ( H^1 ( G , V_{n, d} ) \rightarrow H^1 ( H , V_{n, d} )$ is not trivial, then there exists a non-trivial $[W] \in H^1 ( G, V_{n, d}[p] )$, which is the kernel of the restriction to $H^1 ( H, V_{n, d}[p] )$.
Since $H$ is normal in $G$, we have the inflation-restriction sequence
\[
0 \rightarrow H^1 ( G / H , V_{n, d}[p] ) \rightarrow H^1 ( G , V_{n, d}[p] ) \rightarrow H^1 ( H , V_{n, d}[p] ).
\]
Then $[W]$ is the image by the inflation of a non-trivial element of $H^1 ( G / H , V_{n, d}[p] )$.
Since $H^1 ( G/ H , V_{n, d}[p] ) = 0$, we get a contradiction.
Then $H^1 ( G , V_{n, d} ) \rightarrow H^1 ( H , V_{n, d} )$ is injective.
\CVD

\begin{thm}[Theorem \ref{teo2}]
Suppose that $ G_1 $  contains an element $ g $ whose order divides $ p-1 $ and not fixing any non-trivial element of $ \A[p] $. Moreover suppose that $H^1 ( G_1, \A[p] ) = 0$. Then $H^1_{{\rm loc}} ( G_n , \A[p^n] ) = 0$ for every positive integer $n$ and so the local-global divisibility by $p^n$ holds for $\A ( k )$.
\end{thm}

\DIM Let $n \in \N$ and consider $H^1_{{\rm loc}} ( G_n , \A[p^n] )$.
Let $\widetilde{g} \in G_n$ be such that the restriction of $\widetilde{g}$ to $K_1$ is $g$. 
By applying Corollary \ref{cor12} with $\widetilde{g}$ in the place of $g$, ${\rm Gal} ( K_n / K_1 )$ in the place of $H$ and $\langle \widetilde{g}, H \rangle$ in the place of $G$, we get that $H^1_{{\rm loc}} ( \langle \widetilde{g} , H \rangle , \A[p^n] ) = 0$.
Then for any $[Z] \in H^1_{{\rm loc}} ( G_n , \A[p^n] ) = 0$, $[Z]$ is in the kernel of the restriction to $H^1 ( H , \A[p^n] )$.
Concludes the proof by applying Lemma \ref{lem21}. 
\CVD

\begin{obs}\label{rem22}
We would like to remove the hypothesis on the triviality of $H^1 ( G_1, \A[p] )$ in Theorem \ref{teo2}.
Observe that to do that, by Corollary \ref{cor12} and Remark \ref{rem11}, it would be sufficient to prove the following fact: let $p$ be a prime, let $d$ be a positive integer and $G$ be a subgroup of ${\rm GL}_{2d} ( \Z / p \Z )$, then there exists a $p$-Sylow subgroup of $G$ such that $g$ is in its normalizer.
\end{obs}

In \cite{C-S}, Ciperani and Stix found an interesting relation between the irreducible subquotients of ${\rm End} ( \A[p] )$ and $\A[p]$ as Galois modules and the triviality of a certain Tate-Shafarevich group (see \cite[Theorem 4]{C-S} and the Appendix for the details). 
To study the local-global divisibility problem we need a similar result in which we replace the group studied by Ciperani and Stix with the first local cohomology group.
We do this in the following proposition, that it is also inspired by Section 6 of \cite{LW}.

\begin{pro}\label{pro23}
Let $G$ be a subgroup of ${\rm GL}_{2d} ( \Z / p^n \Z )$ acting on $V_{n, d}$ in the classical way. Let $H$ be the normal subgroup of $G$ of the elements acting like the identity on $V_{n, d}[p]$. 
Suppose that $G$ contains an element $\delta$ such that $\delta -Id$ is a bijective automorphism of $V_{n, d}$ and let $\overline{\delta}$ be its class in $G / H$. Then if $H^1 ( G/ H , V_{n, d}[p] ) = 0$, and  $V_{n, d}[p]$ and ${\rm End} ( V_{n, d}[p] )$ have no common irreducible sub$\Z / p \Z[\langle \overline{\delta} \rangle]$-modules (the action of $\overline{\delta}$ over ${\rm End} ( V_{n, d}[p] )$ is induced by the conjugation), then $H^1 ( G , V_{n, d} ) = 0$.
\end{pro}

\DIM Consider the inflation-restriction sequence
\[
0 \rightarrow H^1 ( G/ H , V_{n, d}[p] ) \rightarrow H^1 ( G, V_{n, d}[p] ) \rightarrow H^1 ( H , V_{n, d}[p] )^{G/ H}.
\]
Since $H^1 ( G / H , V_{n, d}[p] ) = 0$, $H^1 ( G, V_{n, d}[p] )$ is isomorphic to a subgroup of $H^1 ( H , V_{n, d}[p] )^{G/ H}$.
Let $\phi ( H )$ be the Frattini sub-group of $H$ (see the previous subsection for the definition of Frattini sub-group).
In particular recall that $H / \phi ( H )$ is an elementary $p$-abelian group.
Since $H$ acts like the identity over $V_{n, d}[p]$ and $V_{n, d}[p]$ is a commutative group with exponent $p$, $H^1 ( H , V_{n, d}[p] )^{G/ H} = {\rm Hom}_{\Z / p \Z[G/ H]} \{ H/ \phi ( H ) , V_{n, d}[p] \}$ where $G/ H$ has an action induced by coniugacy over $H/ \phi ( H )$.
We shall prove that ${\rm Hom}_{\Z / p \Z[\overline{\delta}]} \{ H/ \phi ( H ) , V_{n, d}[p] \}$, and so ${\rm Hom}_{\Z / p \Z[G/ H]} \{ H/ \phi ( H ) , V_{n, d}[p] \}$, is trivial.

By eventually replacing $\overline{\delta}$ with its $p$-power, we can suppose that $p$ does not divide the order of $\overline{\delta}$.  
Then the action of $\langle \overline{\delta} \rangle$ is semisimple and  $H/ \phi ( H )$ is isomorphic to a direct sum of irreducible $\langle \overline{\delta} \rangle$-modules. 

Take $W$ an irreducible $\Z / p \Z[\langle \overline{\delta} \rangle]$-submodule of $H / \phi ( H )$.
For every non-zero $w \in W$, let $i_w \in \N$ be the biggest integer such that there exists $h \in H$ such that $h \phi ( H ) = w$ and $h \equiv Id \mod ( p^{i_w} )$. 
Then $h \not \equiv Id \mod ( p^{i_w + 1} )$.
Since $W$ is irreducible, every non-zero element of $W$ is a generator of $W$.
Then observe that $i_w$ is the same for every $w \neq 0$.
Thus $W$ is isomorphic to a sub$\Z / p \Z[\langle \overline{\delta} \rangle]$-module of 
\[
M_{i_w + 1} = \{ Id +p^{i_w +1} M \ \mid M \in {\rm Mat}_{2 d} ( \Z / p^n \Z ) \} / \{ Id + p^{i_w +2} M^\prime \ \mid M^\prime \in {\rm Mat}_{2 d} ( \Z / p^n \Z ) \}. 
\]
Since $M_{i_w + 1}$ is isomorphic, as $\Z / p \Z[\langle \overline{\delta} \rangle]$-module, to ${\rm End} ( V_{n, d}[p] )$, $W$ is isomorphic to a submodule of ${\rm End} ( V_{n, d}[p] )$.
Since, by hypothesis, $V_{n, d}[p]$ and ${\rm End} ( V_{n, d}[p] )$ have no common irreducible sub$\Z / p \Z[\langle \overline{\delta} \rangle]$-modules, then ${\rm Hom}_{\Z / p \Z[\overline{\delta}]} \{ H/ \phi ( H ) , V_{n, d}[p] \} = 0$, and so ${\rm Hom}_{\Z / p \Z[G/ H]} \{ H/ \phi ( H ) , V_{n, d}[p] \} = 0$.
Thus $H^1 ( G, V_{n, d}[p] ) = 0$.
Since the groups $H^1 ( G , V_{n, d}[p] )$ and $H^1 ( G, V_{n, d} )[p]$ are isomorphic (see Lemma \ref{obs51}), we get $H^1 ( G, V_{n, d} )[p] = 0$, which implies $H^1 ( G , V_{n, d} ) = 0$.
\CVD  

The following lemma gives a useful criterion to see if an element $\delta$ of $G$ satisfies the hypothesis of Proposition \ref{pro23}.

\begin{lem}\label{lem24}
Let $\delta \in {\rm GL}_{2d} ( \Z / p \Z )$ with order not divisible by $p$ and let $\lambda_1 , \lambda_2 , \ldots , \lambda_{2d}$ the eigenvalues of $\delta$. Suppose that for every $i, j$ between $1$ and $2d$, $\lambda_i / \lambda_j$ is not an eigenvalue of $\delta$. Then $V_{1, d}$ and ${\rm End} ( V_{1, d} )$ have no common irreducible sub$\Z / p \Z[\langle \delta \rangle]$-modules.
\end{lem}

\DIM Observe that the lemma is evident if $\delta$ is diagonalizable over $\F_p$. 
Since $p$ does not divide the order of $\delta$, $\delta$ is diagonalizable in a finite extension $\F_q$ of $\F_p$.
Since the irreducible $\Z / p \Z[\langle \delta \rangle]$-modules are direct sums of irreducible $\F_q[\langle \delta \rangle]$-modules, the result follows.
\CVD

\section{The group of the symplectic similitudes and proof of Theorem 3}\label{sec12}

We start by a description of the Galois action over the $p$-torsion of a polarized abelian variety $\A$ of dimension $d \in \N$.
The referencies that we use for that are \cite[Section 2]{Lom} and \cite{Die}.  

Let $\A$ be an abelian variety admitting a polarization with degree not divisible by $p$. 
The Tate module $T_p ( \A )$ has a skew-symmetric, bilinear, Galois-equivariant form (called Weil pairing)
\[
\langle \ , \ \rangle \colon T_p ( \A ) \times T_p ( \A ) \rightarrow \Z_p ( 1 ),
\] 
where $\Z_p ( 1 )$ is the $1$-dimensional Galois module, in which the action is given by the cyclotomic character $\chi_p \colon {\rm Gal} ( \overline{k} / k ) \rightarrow \Z_p^\ast$.
This is not degenerate over $\A[p]$ because $p$ does not divide the degree of the polarization. 
The fact that the Weil pairing is not degenerate means that the Galois group over $k$ of the field generated by all the torsion points of order a power of $p$ is a subgroup of the group of the symplectic similitude of $T_p ( \A )$ with respect to the Weil pairing ${\rm GSp} ( T_p ( \A ) , \langle \ , \ \rangle )$.
Choosing a basis of $\A[p]$ we can consider $G_1$ (recall that $G_1 = {\rm Gal} ( k ( \A[p] ) / k )$) as a subgroup of ${\rm GSp}_{2d} ( \F_p )$.

For every $\sigma \in G_1$, we define the multiplier of $\sigma$ as the element $\nu ( \sigma ) \in \F_p^\ast$ such that for every $P_1$, $P_2$ in $\A[p]$, $\langle \sigma ( P_1 ) , \sigma ( P_2 ) \rangle = \nu ( \sigma ) \langle P_1 , P_2 \rangle$. 
Then $\nu ( \sigma ) = \chi_p ( \sigma )$ and the determinant of $\sigma$ is just $\nu ( \sigma )^d = \chi_p ( \sigma )^d$.

\begin{thm}[Theorem 3]  
Let $ \A $ be a polarized abelian variety of dimension $ d $ defined over $ k $ and let $ p $ be a prime not dividing the degree of the polarization. 
Suppose that $k \cap \Q ( \zeta_p ) = \Q$.
Set $i = ( ( 2d )! , p-1 )$ and $ k_i $ the subfield of $k ( \zeta_p )$ of degree $i$ over $ k $.
If for every non-zero $P \in \A[p]$ the field $k ( P ) \cap k ( \zeta_p )$ strictly contains $k_i$, then for every positive integer $n$, the group $H^1_{{\rm loc}} ( G_n , \A[p^n] ) = 0$ and so the local-global divisibility by $p^n$ holds for $\A ( k )$.
\end{thm}

\DIM Since $\A$ is a polarized abelian varieties and $p$ does not divide the degree of the polarization, $k ( \zeta_p ) \subseteq K_1$.
Moreover since by hypothesis $k \cap \Q ( \zeta_p ) = \Q$, we have that $\Gal ( k ( \zeta_p ) / k )$ is isomorphic to $\Z / ( p-1 ) \Z$.
Let $N$ be the group $\Gal ( K_1 / k ( \zeta_p ) )$.
By elementary Galois theory, then $N$ is a normal subgroup of $G_1$, containing all $p$-Sylow subgroups of $G_1$ because $[G_1: N] = p-1$, which is not divisible by $p$.
Let $H$ be a $p$-Sylow subgroup of $G_1$.
Let $i$ be $( ( 2d )! , p-1 )$.      
By Corollary \ref{cor14}, there exists $g \in G_1$ of order $( p - 1 )$ such that its restriction to $k ( \zeta_p )$ has order divided by $( p-1 ) / i$. 
By hypothesis, for every point $P$ of order $p$ of $\A$ we have that $k ( P ) \cap k ( \zeta_p )$ strictly contains the subfield of degree $i$ over $k$, which is fixed by the restriction of $g$ to $k ( \zeta_p )$. 
Then $g$ does not fix any point of order $p$ and so $g - Id$ is bijective as endomorphism of $\A[p]$.
Conclude the proof by applying Corollary \ref{cor12}.
\CVD

\section{Proof of Theorem 4}\label{sec3}

The proof of Theorem \ref{teo6} requires the study of some properties of ${\rm GSp}_4 ( \F_p )$.
We do this in the next subsection.

\subsection{Some properties of the group ${\rm GSp}_4 ( \F_p )$}

In the next lemma we list some well-known properties of the group ${\rm GSp}_4 ( \F_p )$. 

\begin{lem}\label{lem141}
Let $p \geq 3$ be a prime number.
\begin{enumerate}
\item  The order of ${\rm GSp}_4 ( \F_p )$ is $p^4 ( p-1 )^3 ( p+1 )^2 ( p^2 + 1 )$;
\item Let $B$ be an element of ${\rm GSp}_4 ( \F_p )$. The eigenvalues of $B$ can be written as $\lambda_1$, $\lambda_2$, $\nu ( B ) \lambda_1^{-1}$, $\nu ( B ) \lambda_2^{-1}$, where $\nu$ is the multiplier (see Section \ref{sec12}).
\end{enumerate}
\end{lem}

\DIM 1. It is well-known.

For 2. we can use that every $M \in {\rm Sp}_4 ( \F_p )$ has eigenvalues $\alpha , \beta , \alpha^{-1} , \beta^{-1}$ (see \cite{Dic}) or \cite[Lemma 2.2]{Die}), and the exact sequence 
\[
1 \rightarrow {\rm Sp}_4 ( \F_p ) \rightarrow {\rm GSp}_4 ( \F_p ) \rightarrow ( \Z / p \Z )^\ast \rightarrow 1,
\]
where the last map is $\nu \colon {\rm GSp}_4 ( \F_p ) \rightarrow ( \Z / p \Z )^\ast$. 
\CVD

The next theorem, proved by Lombardo (see \cite[Section 3.1]{Lom}), gives a very precise list of the maximal subgroups of ${\rm GSp}_4 ( \F_p )$ not containing ${\rm Sp}_4 ( \F_p )$ and it is one of the main ingredient of our proof.

\begin{thm}\label{teo142}
Let $p > 7$ be a prime number. Let $G$ be a proper subgroup of ${\rm GSp}_4 ( \F_p )$ not containing ${\rm Sp_4} ( \F_p )$.
Then $G$ is contained in a maximal proper subgroup $\Gamma$ of ${\rm GSp}_4 ( \F_p ) )$ such that one of the following holds:
\begin{enumerate}
\item $\Gamma$ stabilizes a subspace;
\item There exist $2$-dimensional subspaces $V_1$, $V_2$ of $\F_p^4$ such that $\F_p^4 = V_1 \bigoplus V_2$ and
\[
\Gamma = \{ A \in {\rm GSp}_4 ( \F_p ) \mid \ \exists \gamma \in S_2 \mid \ A V_i \subseteq  V_{\gamma ( i )} \ i = 1, 2 \};
\]
\item There exists a $\F_{p^2}$-structure on $\F_p^4$ such that
\[
\Gamma = \{ A \in {\rm GSp}_4 ( \F_p ) \mid \exists \rho \in {\rm Gal} ( \F_{p^2} / \F_p ) \mid \forall \lambda \in \F_{p^2}, \forall v \in \F_p^4 \ A ( \lambda \ast v ) = \rho ( \lambda ) \ast A ( v ) \},
\]
where $\ast$ is the multiplication map $\F_{p^2} \times \F_p^4 \rightarrow \F_p^4$. In this case, the set 
\[
\{ A \in {\rm GSp}_4 ( \F_p ) \mid \forall \lambda \in \F_{p^2}, \forall v \in \F_p^4 \ A ( \lambda \ast v ) = \lambda \ast A ( v ) \},     
\]
is a subgroup of $\Gamma$ of index $2$;
\item $\Gamma$ contains a group $H$ isomorphic to ${\rm GL}_2 ( \F_p )$ such that the projective image of $\Gamma$ is identical to the projective image of $H$. Moreover for every $\sigma \in H$, the eigenvalues of $\sigma$ can be written as $\lambda_1^3 , \lambda_1^2 \lambda_2 , \lambda_1 \lambda_2^2 , \lambda_2^3$, with $\lambda_1$ and $\lambda_2$ roots of a second degree polynomial with coefficients in $\F_p$. Here $\lambda_1$ and $\lambda_2$ are the eigenvalues of the element of ${\rm GL}_2 ( \F_p )$ corresponding to $\sigma$;
\item The projective image of $\Gamma$ has order at most $3840$.
\end{enumerate}
\end{thm}

\DIM See \cite[Definition 3.1, Definition 3.2, Theorem 3.3, Lemma 3.4]{Lom}.
\CVD

\subsection{Subgroups of ${\rm PGL}_2 ( \F_q )$ and ${\rm SL}_2 ( \F_q )$} 

Let $q$ be a power of $p$. 
To prove Theorem \ref{teo6}, in many cases we can reduce to study a group isomorphic to a subgroup of ${\rm PGL}_2 ( \F_q )$ (see the next subsection), or to a subgroup of ${\rm SL}_2 ( \F_q )$.
Then we recall the well-known classifications of subgroups of  ${\rm PGL}_2 ( \F_q )$ and ${\rm SL}_2 ( \F_q )$ that we often use in the next subsection.

\begin{pro}\label{pro31}   
Let $G$ be a subgroup of ${\rm PGL}_2 ( \F_q )$ of order not divided by $p$. Then, if $G$ is neither cyclic nor dihedral, $G$ is isomorphic to either $A_4$, or $S_4$ or $A_5$.
\end{pro}

\DIM See \cite[Proposition 16]{Ser}.
\CVD

\begin{pro}\label{pro32}
Let $G$ be a subgroup of ${\rm SL}_2 ( \F_q )$ and suppose $p \geq 5$ and $p$ divides the order of $G$. Then either there exists $r \geq 1$ such that $G$ contains ${\rm SL}_2 ( \F_{p^r} )$ or $G$ has a unique abelian $p$-Sylow subgroup $H$ such that $G/ H$ is cyclic of order dividing $q-1$.
\end{pro}

\DIM See \cite[Chapter 3.6, Theorem 6.17]{Suz}.
\CVD

The following corollary of Propositions \ref{pro31} and \ref{pro32} will be often used in the next subsection.

\begin{cor}\label{cor33}
Let $p \geq 5$ be a prime number and let $G$ be a subgroup of ${\rm GL_2} ( \F_p )$ such that $G$ contains an element $\sigma$ of order $> 2$ and dividing $p + 1$, and such that the image of the determinant of $G$ in $\F_p^\ast$ has order $i$. Then $G$ contains a scalar matrix of order $i / ( i, 60 )$.
\end{cor}

\DIM Suppose first that $p$ divides the order of $G$.
Since by hypothesis $\sigma$ has order not dividing $p-1$ and not divided by $p$, by Proposition \ref{pro32} $G$ contains all ${\rm SL}_2 ( \F_p )$.
Since the image of the determinant of $G$ in $\F_p^\ast$ has order $i$, then $G$ contains a scalar matrix of order at least $i / ( i, 2 )$.

Suppose that $p$ does not divide the order of $G$.
Let $\delta \in G$ be such that its determinant has order $i$.
Then, since $i$ divides $p-1$ and $( p-1 , p+1 ) = 2$, by eventually considering a suitable power $g$ of $\delta$, we can suppose that $g$ is diagonalizable and it has determinant of order divided by $i / ( i, 2 )$.
Denote $PG$ the image of $G$ by the projection over ${\rm PGL}_2 ( \F_p )$ and $\overline{g}$, respectively $\overline{\sigma}$, the images of $g$ respectively $\sigma$ in $PG$.
By Proposition \ref{pro31} either $PG$ is cyclic, or $PG$ is dihedral or $PG$ is a group with exponent dividing $60$.

Suppose that $PG$ is cyclic.
Then $\overline{g}$ and $\overline{\sigma}$ commute.
Hence $g \sigma g^{-1} \sigma^{-1}$ is a scalar matrix with determinant $1$.
Since $g$ is diagonalizable and $\sigma$ is not diagonalizable because its order does not divide $p-1$, a simple calculation shows that $g^2$ is a scalar matrix. 
Then $G$ contains a scalar matrix of order $i / ( i, 4 )$.

Suppose that $PG$ is dihedral.
We call a rotation a power of the element of biggest order of $PG$ and a simmetry any element of order $2$ that anticommutes with the rotations.
If $\overline{g}$ and $\overline{\sigma}$ commute, then like in the previous case we prove that $g^2$ is a scalar matrix. 
Moreover, if $\overline{g}$ is a simmetry, then it has order $2$ and so $g^2$ is a scalar matrix.
Then $G$ contains a scalar matrix of order $i / ( i, 4 )$.
Thus it only remains the case where $\overline{\sigma}$ is a simmetry and $\overline{g}$ is a rotation. 
In this case $\overline{\sigma} \overline{g} \overline{\sigma}^{-1} = \overline{g}^{-1}$ and so $\sigma g \sigma^{-1} g$ is a scalar matrix $\mu Id$ with $\mu \in \F_p^\ast$.
Observe that the determinant of $\mu Id$ is equal to the square of the determinant of $g$.
Then also in this case $G$ contains a scalar matrix of order $i / ( i, 4 )$.

It is well-knwown that $A_4$ has exponent $6$, $S_4$ has exponent $12$ and $A_5$ has exponent $30$.
Since $( 30, 12 ) = 60$, and $4$ divides $60$, in particular $G$ contains a scalar matrix of order $i / ( i, 60 )$.
\CVD    

\subsection{End of the proof}

We first recall the statement of Theorem \ref{teo6}.

\begin{thm}[Theorem 4] 
Let $ \A $ be a polarized abelian surface defined over $ k $.
For every prime number $p > 3840$ such that $k \cap \Q ( \zeta_p ) = k$ and not dividing the degree of the polarization, if there exists $n \in \N$ such that $H^1_{{\rm loc}} ( G_n , \A[p^n] ) \neq 0$, then there exists a finite extension $ \widetilde{k} $ of $ k $ of degree $\leq 24$ such that $ \A $ is $\widetilde{k}$-isogenous to an abelian surface with a torsion point of order $ p $ defined over $ \widetilde{k} $.  
\end{thm}

\DIM Suppose that there exists $n \in \N$ such that $H^1_{{\rm loc}} ( G_n , \A[p^n] ) \neq 0$.
The proof is divided in some distinct steps. 
The first is the following simple lemma.

\begin{lem}\label{lem33}
The group $G_1$ is isomorphic to its projective image to ${\rm PGL}_4 ( \F_p )$. Moreover the function $\nu$ from $G_1$ to $( \Z / p \Z )^\ast$ sending $\sigma \in G_1$ to its multiplier $\nu ( \sigma )$ is surjective and $G_1$ contains an element $g$ of order $p-1$ and multiplier divided by $( p-1 ) / 2$.
\end{lem}

\DIM If $G_1$ is not isomorphic to its projective image, then it contains a scalar matrix whose eigenvalue  is distinct of $1$.
Then, by Theorem \ref{teo2}, $H^1_{{\rm loc}} ( G_n , \A[p^n] ) = 0$ for every positive integer $n$ (actually $H^1 ( G_n , \A[p^n] ) = 0$ for every positive integer $n$, see for instance \cite[p. 29]{DZ}).  

Since by hypothesis $k \cap \Q ( \zeta_p ) = \Q$ and $k ( \zeta_p )$ is the subfield of $K_1$ fixed by the kernel of the multiplier $\nu$ (see Section \ref{sec12}), we have $G_1 / \ker ( \nu )$ isomorphic to $\Gal ( k ( \zeta_p ) / k )$ isomorphic to $( \Z / p \Z )^\ast$.

Finally, since $G_1 / \ker ( \nu )$ is a cyclic group of order $( p - 1 )$, every element of $G_1$ whose class generates $G_1 / \ker ( \nu )$ has order divided by $( p-1 )$.
\CVD 

Next proposition shows that a large subgroup of $G_1$ has a stable proper subspace of $\A[p]$.

\begin{pro}\label{pro34}
There exists a subgroup $\Gamma$ of $G_1$ of index at most $4$ and a proper subspace $V$ of $\A[p]$ such that $\sigma ( V ) = V$ for every $\sigma \in \Gamma$.
\end{pro} 

\DIM By Lemma \ref{lem33}, $G_1$ is isomorphic to its projective image and so it does not contain $Sp_4 ( \F_p )$, because $-Id \in Sp_4 ( \F_p )$.
Moreover, see Lemma \ref{lem33}, $G_1$ has order at least $p-1$ and recall that $p > 3840$.
Then, by Theorem \ref{teo142}, either $G_1$ stabilizes a proper subspace of $\F_p^4$, or $G_1$ is contained in a maximal subgroup of type 2., 3., or 4. in the list of Theorem \ref{teo142}.

Suppose that $G_1$ is contained in a subgroup of type 2. Then there exists $V_1$ and $V_2$ subspaces of $\A[p]$ of dimension $2$ such that for every $\sigma \in G_1$, either $\sigma$ permutes $V_1$ and $V_2$ or $\sigma$ stabilizes $V_1$ and $V_2$.
Let $\Gamma$ be the subgroup of $G_1$ that stabilizes $V_1$ and $V_2$.
Observe that it is a normal subgroup of index at most $2$.
Then $\Gamma$ stabilizes two proper subspaces.

Suppose that $G_1$ is contained in a subgroup of type 3. Then $G_1$ has a subgroup $\Gamma$ of index at most $2$ such that there exists a $\F_{p^2}$-structure on $\F_p^4$ such that $\Gamma$ is contained in the group
\[
\{ A \in {\rm GSp}_4 ( \F_p ) \mid \forall \lambda \in \F_{p^2}, \forall v \in \F_p^4 \ A ( \lambda \ast v ) = \lambda \ast A ( v ) \},
\]
where $\ast$ is the multiplication map $\F_{p^2} \times \F_p^4 \rightarrow \F_p^4$.
Then, by choosing a $\F_p^2$-basis of $\F_p^4$, we get an injective homomorphism of $\phi \colon \Gamma \rightarrow {\rm GL}_2 ( \F_{p^2} )$.
Also observe that for every $\sigma \in \Gamma$, $\phi ( \sigma )$ has the same eigenvalues of $\sigma$ (with multiplicity divided by $2$).
Then $\phi ( \Gamma )$ is contained in ${\rm PGL}_2 ( \F_{p^2} )$.
Suppose first that $p$ does not divide the order of $\Gamma$.
Then by Proposition \ref{pro31} and the fact that $p-1$ divides the order of $G_1$, either $\Gamma$ is cyclic or $\Gamma$ is dihedral. 
If $\Gamma$ is cyclic, then, since the generator of $\Gamma$ has two eigenvalues with multiplicity $2$, it stabilizes two subspaces of dimension $2$.
If $\Gamma$ is dihedral, then it contains a normal cyclic subgroup $\Gamma^\prime$ of index $2$.
Thus, by replacing $\Gamma$ with $\Gamma^\prime$, we reduce us to the previous case.
Also observe that $[G_1: \Gamma^\prime]$ divides $4$.
Suppose now that $p$ divides the order of $\Gamma$.
Then, by Proposition \ref{pro32} and the fact that $\phi ( \Gamma )$ is isomorphic to its projective image, $\phi ( \Gamma ) \cap {\rm SL}_2 ( \F_{p^2} )$ is contained in a Borel subgroup.
Since the $p$-Sylow is normal, actually $\phi ( \Gamma )$ is contained in a Borel subgroup and so $\Gamma$ stabilizes a subspace of dimension $2$.

Suppose that $G_1$ is contained in a maximal subgroup of type 4. 
Then, since $G_1$ is isomorphic to its projective image, $G_1$ is isomorphic to a subgroup of ${\rm GL}_2 ( \F_p )$.
Observe that (see Theorem \ref{teo142}) the isomorphism sends the projective image of $G_1$ to ${\rm PGL}_2 ( \F_p )$ and so actually $G_1$ is isomorphic to a subgroup of ${\rm PGL}_2 ( \F_p )$.
If $p$ does not divide the order of $G_1$, then by Proposition \ref{pro31} and since $p-1$ divides the order of $G_1$, we get that $G_1$ is either cyclic or dihedral. 
If $G_1$ is cyclic and since ${\rm PGL}_2 ( \F_p )$ has order $p ( p- 1 ) ( p+1 )$, we get that $G_1$ stabilizes a subspace.
If $G_1$ is dihedral, then $G_1$ has a normal cyclic subgroup $\Gamma$ of index $2$ and so, by replacing $G_1$ with $\Gamma$, we get the same result.
Suppose that $p$ divides the order of $G_1$.
In this case, by Proposion \ref{pro32}, $G_1$ has a unique non-trivial $p$-Sylow subgroup and so it stabilizes a subspace.
\CVD

From the next proposition and a deep result of Katz (see Theorem \ref{teogrande}) it will easily follow Theorem \ref{teo6}.  

\begin{pro}\label{pro35}
There exists a subgroup $\Gamma$ of $G_1$ of index $\leq 24$ such that every $\gamma \in \Gamma$ has at least an eigenvalue equal to $1$.
\end{pro}

\DIM By Proposition \ref{pro34}, by eventually replacing $G_1$ with a subgroup $\Gamma$ of index $2$ or $4$, there exists $V$ a proper subspace of $\A[p]$ stable by the action of $\Gamma$. 
Then, by Lemma \ref{lem33}, $\Gamma$ contains a diagonal element $g$ with order dividing $p-1$ and multiplier with order divided by $( p-1 ) / ( p-1 , 8 )$. 
By abuse of notation, from now on we set $G_1 = \Gamma$.
Let $V^{\perp}$ be the subspace of $\A[p]$ of the $w \in \A[p]$ such that, for every $v \in V$, we have $\langle v , w \rangle = 0$.
Then, since the Weil pairing $\langle \ , \ \rangle$ is Galois-equivariant, also $V^{\perp}$ is stable by the action of $G_1$.
Suppose first that $V$ has dimension $1$.
Then $V \subseteq V^{\perp}$ and $V^{\perp}$ has dimension $3$.
On the other hand, if $V$ has dimension $3$, then $V^\perp$ has dimension $1$ and $V^\perp \subseteq V$.
By eventually replacing $V$ with $V^\perp$, we have two cases: either $V$ has dimension $3$ or $V$ has dimension $2$.
\vspace{5 pt}

{\bf The case when $V$ has dimension $3$.} Suppose that $V$ has dimension $3$ and so $V^{\perp}$ has dimension $1$ and it is contained in $V$.
Then we have the following $G_1$-modules: $V^{\perp} \subseteq V \subseteq \A[p]$ of dimension $1$, $V \subseteq \A[p]$ of dimension $3$, $V / V^{\perp}$ of dimension $2$ and $\A[p] / V$ of dimension $1$.
In particular observe that the exponent of $G_1$ is coprime with $( p^2 + 1 )/2$.
Let $H$ be a $p$-Sylow of $G_1$.
Then $H$ is the identity over $V^{\perp}$ and over $\A[p]/ V$.
Then for every $\tau \in G_1$, if the projection of $\tau$ over $V / V^{\perp}$ is in the normalizer of the projection $H$, then $\tau$ is in the normalizer of $H$.
Since $V / V^{\perp}$ has dimension $2$ and for every subgroup $\Delta$ of ${\rm GL}_2 ( \F_p )$, every element of order dividing $p-1$ is in the normalizer of a $p$-Sylow subgroup of $\Delta$, every element of $G_1$ of order dividing $p-1$ is in the normalizer of a $p$-Sylow subgroup of $G_1$.
Then, see Corollary \ref{cor12} and Remark \ref{rem11}, every element of order dividing $p-1$ has at least an eigenvalue equal to $1$.
Let $\sigma$ be in $G_1$ such that $\sigma$ has all the eigenvalues distinct from $1$.
Then, since $\sigma$ stabilizes $V^{\perp}$ and $A[p] / V$, the unique possibility is that the automorphism of $V / V^{\perp}$ induced by $\sigma$ has order divided by a divisor of $( p+1 )$ not dividing $( p-1 )$.
On the other hand choose $v$ and $w$ in $V$ such that $\{ v , w \}$ is sent by the projection to a basis of $V / V^{\perp}$. 
Let us remark that $v$ is not orthogonal to $w$. 
In fact, if $v$ would not be orthogonal to $w$, then $\langle v \rangle^{\perp}$ would be equal to $V$, and so $V^{\perp}$ would be $\langle v \rangle$.
But $v \not \in V^{\perp}$.
Then we have a contradiction.
Thus for every $\tau \in G_1$, the determinant of the projection of $\tau$ over $V / V^{\perp}$ is equal to the multiplier of $\tau$ and so, by Corollary \ref{cor33} and the fact that the image of the multiplier has index dividing $4$ over $\F_p^\ast$, we have that there exists $\delta \in G_1$ such that the projection of $\delta$ over $V / V^{\perp}$ is a scalar matrix $\lambda Id$ with $\lambda$ of order $( p-1 ) / ( p-1 , 60 )$.
Then, since the order of $\delta$ divides $p-1$, one of its eigenvalue is $1$. 
Then the eigenvalues of $\delta$ are $1, \lambda, \lambda, \lambda^2$.
Observe that the eigenvalues distinct from $\lambda$ are the eigenvalue of the restriction of $\delta$ to $V^{\perp}$ and the other to the projection of $\delta$ to $\A[p] / V$.
Suppose that $\delta$ is the identity over $V^{\perp}$ (the other case is identical).   
Let $\gamma \in G_1$ be any element of order dividing $p-1$ and suppose that the eigenvalues of $\gamma$ are $\lambda_1 , \lambda_2 , \lambda_3 , \lambda_4$. 
Then observe that since the projection of $\delta$ to $V / V^{\perp}$ is in the center of the projection of $G_1$, by eventually permuting the eigenvalues of $\gamma$, for every integer $i$ we have that the eigenvalues of $\delta^i g$ are $\lambda_1 , \lambda^i \lambda_2 , \lambda^i \lambda_3 , \lambda^{2i} \lambda_4$.
Moreover $\delta^i \gamma$ has order dividing $( p-1 ) p^r$ for a certain integer $r$.
But raising a power of $p$ of an element does not change the eigenvalues and so we can suppose that $\delta^i \gamma$ has order dividing $p-1$.
Since $\lambda$ has order $p-1$ and $p > 3840$, if $\lambda_1 \neq 1$ then we can choose $i$ such that $\delta^i \gamma$ has all eigenvalues distinct from $1$.
But this is not possible and so every element of order dividing $p-1$ is the identity over $V^{\perp}$.
Let again $\sigma$ be an element with all eigenvalues distinct from $1$ and so such that $\overline{\sigma}$ has order divided by a divisor of $( p+1 )$ not dividing $( p-1 )$.
Since we can suppose that $p$ does not divide the order of $\sigma$, then $\sigma^{p+1}$ has order dividing $( p-1 )$.
Thus $\sigma^{p+1}$ is the identity over $V^{\perp}$.
But $( p+1 , p-1 ) = 2$ and so the restriction of $\sigma$ to $V^{\perp}$ is either the identity or $-Id$.
Thus the subgroup $\Gamma$ of $G_1$ that fixes $V^{\perp}$ has index $2$ and so this concludes the proof in the case that $V$ has dimention $3$.       
\vspace{5 pt}

{\bf The case when $V$ has dimension $2$ and $V \cap V^{\perp} = \{ 0 \}$.} Since $V \cap V^{\perp} = \{0\}$, then $\A[p]$ is isomorphic as $G_1$-module to the direct sum of $V$ and $V^{\perp}$.
Moreover we can suppose that $V$ and $V^{\perp}$ are irreducible because, if not, $\A[p]$ has $G_1$-submodule of dimension $1$ and we are in the previous case. 
Suppose that the order of $G_1$ is coprime with $( p+1 ) / 2$.
Then $G_1$ has a unique $p$-Sylow and, by Corollary \ref{cor7} and Remark \ref{rem11}, we have that every element of $G_1$ of order dividing $p-1$ has at least an eigenvalue equal to $1$.
Since $G_1$ has order coprime with $( p+1 ) / 2$ and it stabilizes two spaces of dimension $2$, $G_1$ has exponent dividing $( p-1 ) p^2$.
Since for every $\tau \in G_1$, $\tau$ and $\tau^p$ have the saime eigenvalues, all elements of $G_1$ has at least an eigenvalue equal to $1$.
Then we can suppose that there exists $\sigma \in G_1$ of order dividing $p+1$ and not dividing $( p-1 )$.
 In particular the restriction of $\sigma$ to either $V$ or $V^{\perp}$ should have the same property and so suppose that this is the case for the restriction to $V$ (the other case is identical).
Since $V \cap V^{\perp} = \{ 0 \}$, for every $\tau \in G_1$ the determinant of the restriction of $\tau$ to $V$ is the multiplier of $\tau$.
Since the multiplier has index dividing $4$ over $\F_p^\ast$, by Corollary \ref{cor33} there exists $\delta \in G_1$ such that the restriction of $\delta$ to $V$ is a scalar matrix $\lambda Id$ with $\lambda$ of order $( p-1 ) / ( p-1 , 60 )$.
By eventually replacing $\delta$ with its power, since $( p-1 , p+1 ) = 2$, we can suppose that $\delta$ has order dividing $p-1$, but then we have just that $\lambda$ has order divided  by $( p-1 ) / ( p-1 , 120 )$   
.
In particular observe that since the restriction of $\delta$ to $V$ is a scalar matrix, $\delta$ is diagonalizable over $V^{\perp}$ and $V^{\perp}$ has dimension $2$, then $\delta$ is in the normalizer of a $p$-Sylow subgroup of $G_1$.
Hence, by Corollary \ref{cor12} and Remark \ref{rem11}, the eigenvalues of the restriction of $\delta$ to $V^{\perp}$ are $1$ and $\lambda^2$.
Consider now the restriction $G_{1, \perp}$ of $G_1$ to $V^{\perp}$.
If there exists $\tau \in G_1$ whose restriction to $V^{\perp}$ has order dividing $p+1$ and not divided by $p-1$, then by Corollary \ref{cor33} there exists $\delta^\prime \in G_1$, which is a scalar matrix over $V^{\perp}$ with order dividing $( p-1 ) / ( p-1 , 120 )$.
Then $\delta^\prime$ commutes with $\delta$ and, since $p > 3840$, by making the product of suitable powers of $\delta$ and $\delta^\prime$ we get an element of $G_1$ of order dividing $p-1$ with all the eigenvalues distinct from $1$. 
Moreover, since $V$ and $V^{\perp}$ have dimension $2$, such an element is in the normalizer of a $p$-Sylow subgroup and so $H^1_{{\rm loc}} ( G_n , \A[p^n] ) = 0$ for every positive integer $n$, by Corollary \ref{cor12} and Remark \ref{rem11}.
Then the restriction of $G_{1 , \perp}$ has order dividing $2 p ( p-1 )^2$.
If $p$ divides the order of $G_{1 , \perp}$, by Proposition \ref{pro32} either $G_{1 , \perp}$ contains ${\rm SL}_2 ( \F_p )$ (and so $p+1$ divides the order of $G_{1, \perp}$) or $G_{1, \perp}$ has a unique $p$-Sylow of order $p$.
But in the last case $V^{\perp}$ is reducible and then we get a contradiction.
Hence $G_{1, \perp}$ has order dividing $2 ( p-1 )^2$.
Since $V^{\perp}$ is irreducible, the unique possibility is that $G_{1, \perp}$ has a commutative normal subgroup $\Delta$ of index $2$ with order dividing $( p-1 )^2$.
Take $\Gamma$ the subgroup of $G_1$ of the elements whose restriction to $V^{\perp}$ is in $\Delta$.
Then $G_1$ has index $2$ over $\Gamma$ and $\Gamma$ stabilizes a subspace of dimension $1$ and its perpendicular (then a space of dimension $3$).
Then we are in the previous case already studied: the case when $V$ has dimension $3$.
\vspace{5 pt}

{\bf The case when $V = V^{\perp}$}. First observe that if $V$ is not irreducible, then we are in the case when $V$ has dimension $3$ (or $1$) and so we suppose that $V$ is an irreducible $G_1$-module.
Let $W$ be the $G_1$-module $\A[p] / V$.
Let us call $I_V$, respectively $I_W$, the normal subgroups of $G_1$ fixing all the elements of $V$, respectively $W$.
Suppose that $p$ divides the order of $G_1 / I_V$.
Then there is $\sigma \in G_1$ of order $p$ and a basis $\{ v_1 , v_2 \}$ of $V$ such that $\sigma ( v_1 ) = v_1$ and $\sigma ( v_2 ) = v_1 + v_2$.
Let $w_1$, $w_2$ be in $\A[p]$, such that $w_i$ is not orthogonal to $v_i$ and $w_i$ is orthogonal to $v_j$ for $i \neq j$ and $i, j \in \{1, 2 \}$.
Then $\{v_1 , v_2 , w_1 , w_2 \}$ is a basis of $\A[p]$.
Moreover let $\overline{w_1}$ and $\overline{w_2}$ be the class of $w_1$, respectively $w_2$ modulo $V$.
Then $\{ \overline{w_1} , \overline{w_2} \}$ is a basis of $W$.
Let us show that the class of $\sigma$ in $G_1 / I_W$ has order $p$.
If not $\sigma$ should be the identity.
Then it would exist $v \in V$ such that  $\sigma ( w_1 )$ should be equal to $w_1 + v$.
Thus 
\[
\langle \sigma ( v_2 ) , \sigma ( w_1 ) \rangle = \langle v_1 + v_2 , w_1 \rangle = \langle v_1 , w_1 \rangle.
\]
But $\langle \sigma ( v_2 ) , \sigma ( w_1 ) \rangle = \langle v_2 , w_1 \rangle$, which is distinct from   
$\langle v_1 , w_1 \rangle$ because $v_1$ and $w_1$ are not orthogonal and $v_2 , w_1$ are orthogonal.
In the same way we can prove that if $p$ divides $G_1 / I_W$, then there exists $\sigma \in G_1$ of order $p$ such that the restriction of $\sigma$ to $V$ has order $p$.
Since $V$ and $W$ are irreducible, if their $p$-Sylow is not the identity, their $p$-Sylow cannot be normal and so, by Proposition \ref{pro32}, $G_1 / I_V$ and $G_1 / I_W$ contain all the group ${\rm SL}_2 ( \F_p )$.
Then observe that there exists $\tau_1 \in G_1$ whose restriction over $V$ is $-Id$ and $\tau_2 \in G_1$ whose projection over $W$ is $-Id$.
By Lemma \ref{lem141}, then the other eigenvalues of $\tau_1$ are identical and so a $p$ power of $\tau_1$ is a diagonal matrix with two eigenvalues equal to $-1$ and the others equal to a $\lambda \in \F_p^\ast$.
In the same way we can prove that a $p$-power of $\tau_2$ is a diagonal matrix with two eigenvalues equal to $-1$ and the other equal to $\mu$ for a certain $\mu \in \F_p^\ast$.
Then either $\tau_1$, or $\tau_2$ or $\tau_1 \tau_2$ has order dividing $p-1$, it has all the eigenvalues distinct from $1$ and it is in the normalizer of a $p$-Sylow because it is a scalar matrix over $V$ and over $W$.
Then by Corollary \ref{cor12} and Remark \ref{rem11}. for every positive integer $n$ we have $H^1_{{\rm loc}} ( G_n , \A[p^n] ) = 0$.
Then $G_1 / I_V$ and $G_1 / I_W$ have orders not divided by $p$ and so $G_1$ has a normal $p$-Sylow $N$ such that $N \subseteq I_V$ and $N \subseteq I_W$.
Thus if $G_1$ has order coprime with $( p+1 ) /2$, by Corollary \ref{cor12} and Remark \ref{rem11}. for every positive integer $n$ we have $H^1_{{\rm loc}} ( G_n , \A[p^n] ) = 0$.
Then there exists $\sigma \in G_1$ with order dividing $p+1$ and not dividing $p-1$.
Since $( p-1 , p+1 ) = 2$, by Lemma \ref{lem141} we can suppose that the eigenvalues of $\sigma$ are either $\mu$, $\mu^p$ (both with multiplicity $2$) or $\mu$, $\mu^p$, $-\mu$, $-\mu^p$.
The following lemma, whose proof is similar to the proof of Proposition \ref{pro23}, gives a strong restriction to the order of $\sigma$.

\begin{lem}\label{lemsecondo}
If there exists $n \in \N$ such that $H^1_{{\rm loc}} ( G_n , \A[p^n] ) \neq 0$, then $\sigma$ has order dividing $6$.
\end{lem}

\DIM First observe $\sigma - Id$ is bijective as endomorphism of $\A[p]$ because $\mu \not \in \F_p$.
Moreover, by using Lemma \ref{lem24}, we can prove that $\A[p]$ and ${\rm End} ( \A[p] )$ have a common $\Z / p \Z [ \langle \sigma \rangle]$-modulo only if $\sigma$ has order dividing $6$.
Then by Proposition \ref{pro23} if $H^1 ( G_1 , \A[p] ) = 0$, we immediately get the result. 
Then let us prove that $H^1 ( G_1 , \A[p] ) = 0$ (actually the proof is similar to the proof of Proposition \ref{pro23}).
Consider the exact sequence of $G_1$-modules:
\[
0 \rightarrow V \rightarrow \A[p] \rightarrow W \rightarrow 0
\]
where the first map is the inclusion and the second is the projection.
Since $\delta - Id$ is bijective over $\A[p]$, we have $H^0 ( G_1 , W ) = 0$ and so we get the following long cohomology exact sequence
\[
0 \rightarrow H^1 ( G_1 , V ) \rightarrow H^1 ( G_1 , \A[p] ) \rightarrow H^1 ( G_1 , W ).
\]
Then, to prove the triviality of $H^1 ( G_1 , \A[p] )$, it is sufficient to prove the triviality of $H^1 ( G_1 , V )$ and $H^1 ( G_1 , W )$.
Let us prove the triviality of $H^1 ( G_1 , V )$ (the proof of the triviality of $H^1 ( G_1 , W )$ is identical).
Recall that $N$ is the $p$-Sylow of $G_1$ and $N$ fixes $V$ and $W$.
Then we have the following inflation-restriction sequence:
\[
0 \rightarrow H^1 ( G_1 / N , V ) \rightarrow H^1 ( G_1 , V ) \rightarrow H^1 ( N , V )^{G_1 / N}.
\]
Since $p$ does not divide the order of $G_1/ N$, we have $H^1 ( G_1/ N , V ) = 0$.
Since $N$ fixes $V$, we get that $H^1 ( N, V )^{G_1 / N}$ is isomorphic to ${\rm Hom}_{\Z / p \Z[G_1 / N]}    ( N , V )$ where $G_1$ acts over $N$ by conjugation (recall that since $N$ fixes $V$ and $W$, $N$ is an abelian group with exponent dividing $p$).
By Lemma \ref{lem24}, the action of $\delta$ by conjugation over $N$ is given by an automorphism with eigenvalues contained in the set either $\{ 1 , \mu^{p-1} , \mu^{1-p} \}$, or $\{1, -1,  \mu^{p-1}, -\mu^{p-1}, \mu^{1-p}, -\mu^{1-p} \}$.
On the other hand, over $V$ the element $\delta$ has eigenvalues either $\{\mu , \mu^p\}$ or $\{ \mu, -\mu, \mu^p , -\mu^p\}$. 
But $\{1, -1,  \mu^{p-1}, -\mu^{p-1}, \mu^{1-p}, -\mu^{1-p} \} \cap \{\mu , -\mu , \mu^p , -\mu^p \}$ is not empty only if $\mu$ has order dividing $6$.
Then if $\sigma$ has not order dividing $6$, $H^1 ( G_1 , \A[p] ) = 0$.
\CVD

Observe that if $\sigma$ has order $3$ or $6$, then $\sigma^2$ has order $3$ and it has eigenvalues $\lambda$, $\lambda^p$ (both with multiplicity $2$) and $\lambda$ of order $3$.
Now recall that $G_1$ contains an element $g$ of order dividing $p-1$ and multiplier divided by $( p - 1 ) / ( p-1 , 8 )$.
By Corollary \ref{cor12} and Remark \ref{rem11} $G_1$ has at least an eigenvalue equal to $1$.
Suppose that the corresponding eigenvector is in $V$ (the case when it is in $W$ is identical).
By Proposition \ref{pro31}, since $p$ does not divide the order of $G_1 / I_V$, the projective image of $G_1 / I_V$ is either cyclic, or dihedral or isomorphic to an exceptional subgroup (either $A_4$, or $S_4$, or $A_5$).
If this last case is verified, then $G_1 / I_V$ contains an element $\tau$ which act like $-Id$ over $V$.
By Corollary \ref{cor12} and Remark \ref{rem11}, it acts like the identity over $W$.
Then a suitable $p$-power of $\tau$ commutes with $g$ and by choosing $i = 1$ or $2$, $g^i \tau$ has all the eigenvalues distinct from $1$, because $g$ has multiplier divided by $( p-1 ) / ( p-1 , 8 )$ and $p > 3840$.
Thus either the projective image of $G_1 / I_V$ is cyclic of order $3$ or it is dihedral of order $6$ (in the two cases generated by the class of $\sigma^2$ of order $3$ and the class of $g$ that can have order $1$ or $2$).
Since $g$ has an eigenvalue equal to $1$ over $V$ the unique possibility is that $g$ is either the identity or it has order $2$ over $V$. 
Then $G_1 / I_V$ is generated by $\sigma^2$, $g$ and eventually $\delta \in G_1$, which is a scalar matrix over $V$.
But in this case take a suitable power of $g^2$ multiplied by $\delta$ and get a matrix with order divided $p-1$ and all the eigenvalue distinct from $1$.
Then $G_1$ has either index $3$ or index $6$ over $I_V$ and so
$K_1^{I_V}$ is an extension of degree dividing $6$ of $k$ in which all the elements of $I_V$ fixes all the elements of a subspace of $V$ of dimension $2$.
\CVD  

The last result we need to finish the proof is the following deep result of Katz.

\begin{thm}\label{teogrande}
Let $\B$ be an abelian surface defined over a number field $F$.
If for all but finitely many prime numbers $r$, we have that a prime number $q$ divides the order of $\A ( \F_r )$, then there exists $\B^\prime$ abelian surface defined over $F$ and $F$-isogenous to $\B$ such that $\B^\prime$ admits a point of order $q$ defined over $F$.
\end{thm}

\DIM See \cite[Introduction]{Kat}.
\CVD  

Observe that if all elements of $G_1$ have at least an eigenvalue equals to $1$, then by Chebotarev density theorem for all but finitely many prime numbers $q$, we have that $p$ divides the order of $\A ( \F_q )$. 
By Proposition \ref{pro35}, there exists an extension $L$ of $k$ of degree $\leq 24$ such that every element of $\Gal ( K_1 / L )$ fixes at least a non-trivial element of $\A[p]$.
By applying Theorem \ref{teogrande}, we then conlude the proof.
\CVD 

\section{The counterexample}\label{sec6}

Let $p$ be a prime number such that $p \equiv 2 \mod ( 3 )$. 
Consider the following subgroups of ${\rm GL}_2 ( \Z / p^2 \Z )$:
\[
H_2 = 
\bigg \{
h ( a, b ) =
\left(
\begin{array}{cc}
1 + p ( a - 2b ) & 3p ( b-a ) \\
-pb & 1 - p ( a - 2b ) \\
\end{array}
\right)  
\ a, b \in \Z/ p^2 \Z
\bigg \} 
\]
and 
\[
G_2 =
\bigg \langle
g =
\left( 
\begin{array}{cc} 
1 & -3 \\
1 & -2 \\
\end{array}
\right)
, \
H_2
\bigg \rangle.
\]
A simple computation gives that $g$ has order $3$, which does not divide $p-1$.
A simple verification gives that for every $a, b$, we have 
\[
gh ( a, b ) g^{-1} = h ( -b, a-b ), \ g^2h ( a, b ) g^{-2} = h ( b-a, -a ).
\]
Then $H_2$ is a normal abelian subgroup of $G_2$.     

We shall prove that $H^1_{{\rm loc}} ( G_2, ( \Z /p^2 \Z )^2 ) \neq 0$, by explicitly costructing a cocycle from $G_2$ to $( \Z / p^2 \Z )^2$ that satisfies the local conditions, but it is not a coboundary.
Observe that $H_2$ is a $\Z / p \Z [\langle g \rangle]$-module, with $g$ that acts by conjugation. 
Let $Z$ be a cocycle from $G_2$ to $( p \Z/ p^2 \Z )^2$. 
By cocycle relations and the fact that $H_2$ acts in the same way of the identity over $( p \Z/ p^2 \Z )^2$, we have that $Z$ is homomorphism of $\Z / p \Z[\langle g \rangle]$-modules from $H_2$ to $( p \Z/ p^2 \Z )^2$.   
Using that $g h ( a, b ) g^{-1} = h ( -b, a-b )$, a simple computation shows that the group of homomorphisms of $\Z / p \Z [\langle g \rangle]$-modules from $H_2$ to $( p \Z / p^2 \Z )$ is cyclic generated by $Z \colon H_2 \rightarrow ( p \Z / p^2 \Z )^2$ with $Z_{h ( a, b )} = ( p ( a-2b ), p ( a-b ) )$.
Then, extending $Z$ to $G_2$ by sending $g$ to $( 0, 0 )$ and using the properties of cocycles, we have a cocycle from $G_2$ to $( \Z /p^2 \Z )^2$.
  
Let us show that $Z$ satisfies the local conditions.
In other words we shall prove that for every $( a, b )$ the system $h ( a, b ) - Id ( x, y ) = Z_{h ( a, b )}$ has a solution.
Observe that, by definition of $h ( a, b )$, it is sufficient to prove that if $a \neq 0$ or $b \neq 0$, then
\[   
\left(
\begin{array}{cc}
a - 2b & 3 ( b-a ) \\
-b & 2b - a \\
\end{array}
\right)
\]
has determinant distinct from $0$ in $\Z / p \Z$.
A simple computation show that the determinant is $\Delta ( ( a, b ) ) = a^2 + b^2 -ab$.
Since $\Delta ( ( a, b ) )$ is a homogenous polynomial in $a$ and in $b$ and $a$ and $b$ are symmetric, if it has a non-zero solution, it has a solution of the form $( -1 , \beta )$.
Then $\beta^2 + \beta + 1 = 0$ that gives that $\beta$ has order $3$ in $( \Z / p \Z )^\ast$.
This is not possible because $p \equiv 2 \mod ( 3 )$.
Then $Z$ satisfies the local conditions.

We show that $Z$ is not a coboundary.
Observe that $( h ( 1, 1 ) -Id ) ( x, y ) = Z_{h ( 1, 1 )}$ if and only if $( x, y ) = ( 1, 1 )$.
Moreover $( h ( 2, 1 ) -Id ) ( x, y ) = Z_{h ( 2, 1 )}$ if and only if $( x, y ) = ( -1, 0 )$.
Then $Z$ is not a coboundary.

Let $k$ be a number field and let $\E$ be a not $CM$ elliptic curve defined over $k$.
By the main result of \cite{Ser}, for every prime number $l$ sufficiently big the representation of ${\rm Gal} ( \overline{k} / k )$ over the group of the automorphism on the Tate $l$-module of $\E$ is surjective. 
Choose a prime $p \equiv 2 \mod ( 3 )$ sufficiently big.
Then ${\rm Gal} ( k ( \E[p^2] ) / k ) = {\rm GL}_2 ( \Z / p^2 \Z )$.
Let $L$ be the field contained in $k ( \E[p^2] )$ fixed by $G_2$.
Then ${\rm Gal} ( L ( \E[p^2] ) / L )$ is isomorphic to $G_2$ and the isomorphism induces an isomorphism of the ${\rm Gal} ( L ( \E[p^2] ) / L )$-module $\E[p^2]$ over the $G_2$-module $( \Z / p^2 \Z )^2$. 
Then $H^1_{{\rm loc}} ( {\rm Gal} ( L ( \E[p^2] ) / L ) , \E[p^2] )$ is not trivial.
By \cite[Theorem3]{DZ3}, by eventually replacing $L$ with a field $L^\prime$ such that $L \subseteq L^\prime$ and $L ( \E[p^2] ) \cap L^\prime = L$, we get a counterexample on local-global divisibility by $p^2$ over $\E ( L^\prime )$.
Observe that ${\rm Gal} ( L^\prime ( \E[p] ) / L^\prime )$ is generated by an element of order $3$, then of order not dividing $p-1$.
Moreover $H^1 ( {\rm Gal} ( L^\prime ( \E[p] ) / L^\prime ) , \E[p] ) = 0$ because $p$ and $3$ are distinct.     
   
\section{Appendix}\label{sec8}

Ciperani and Stix \cite{C-S} studied the following question of Cassels, which is related to the local-global divisibility problem: let $ k $ be a number field and let $ \A $ be an abelian variety defined over $k$. 
For every prime number $q$ we say that the Tate-Shafarevich group $\Sh ( \A / k )$ is $q$-divisible in $H^1 ( k , \A )$ if $\Sh ( \A / k ) \subseteq \cap_{n \in \N^\ast} q^n H^1 ( k , \A )$. 
What is the set of prime numbers $q$ such that $\Sh ( \A / k )$ is $q$-divisible ?

We explain the criterion found by Ciperani and Stix to answer to this question. 
Define 
\[
\Sh^1 ( k, \A [p^n] ) = \cap_{v \in M_k} \ker ( H^1 ( \Gal ( \overline{k}/ k ), \A[p^n] ) \rightarrow H^1 ( \Gal ( \overline{k}_v/ k_v ), \A[p^n] ) ).
\]
Let $\A^t$ be the dual variety of $\A$. 
Ciperani and Stix (see \cite[Proposition 13]{C-S}) proved the following result. 

\begin{thm}\label{teo311}
If $\Sh^1 ( k, \A^t [p^n] )$ is trivial for every positive integer $n$, then $\Sh ( \A / k )$ is $p$-divisible over $H^1 ( k, \A )$.
\end{thm} 

Then in their paper they found very interesting criterions for the triviality of $\Sh^1 ( k, \A^t [p^n] )$ (see \cite[Theorems A, B, C, D]{C-S}) and we applied some of their ideas in this paper, in particular in the second subsection of Section \ref{sec2}.
In particular they got a very strong criterion in the case when $\A$ is an elliptic curve.

We now explain the relation with Cassels question and the local-global divisibility problem (observe that Ciperani and Stix \cite[Remark 20]{C-S} already substantially observed the connection. Here we want just to precise it).
Let $\Sigma$ be a subset of the set of places $M_k$ of $k$.
By following \cite[p. 15]{San} with $G = \A[p^n]$, we define 
\[
\Sh^1_\Sigma ( k, \A [p^n] ) = \cap_{v \not \in \Sigma} \ker ( H^1 ( \Gal ( \overline{k}/ k ), \A[p^n] ) \rightarrow H^1 ( \Gal ( \overline{k}_v/ k_v ), \A[p^n] ) ),
\]
\[
\Sh^1_\omega ( k, \A [p^n] ) = \cup_{\Sigma \ {\rm finite}} \Sh^1_\Sigma ( k, \A [p^n] ).   
\]
Observe that $\Sh^1 ( k , \A[p^n] ) = \Sh^1_{\emptyset} ( k , \A[p^n] )$ and obviously $\Sh^1 ( k , \A[p^n] ) \subseteq \Sh^1_\omega ( k , \A[p^n] )$.
The Lemma \cite[Lemme 1.2]{San} applied with $B = \A[p^n]$ implies that $\Sh^1_\omega ( k , \A[p^n] )$ is isomorphic to $H^1_{{\rm loc}} ( G_n , \A[p^n] )$.
Then if $H^1_{{\rm loc}} ( G_n , \A[p^n] ) = 0$, we have $\Sh^1 ( k , \A[p^n] ) = 0$.
By Theorem \ref{teo311} we then get the following corollaries of respectively Theorem \ref{teo2}, Theorem \ref{teo5} and Theorem \ref{teo6}. 

\begin{cor}\label{cor3}  
Suppose that $\Gal ( k ( \A^t [p] ) / k )$ contains an element $ g $ whose order divides $p-1$ and not fixing any non-trivial element of $\A^t[p]$. Moreover suppose that $H^1 ( \Gal ( k ( \A^t[p] ) / k ) , \A^t[p] ) = 0$. Then $\Sh ( \A / k )$ is $p$-divisible in $H^1 ( k , \A )$.
\end{cor}  
  
\begin{cor}\label{cor6}
Let $\A$ be a principally polarized abelian variety of dimension $d$ defined over $k$ and suppose that $k \cap \Q ( \zeta_p ) = \Q$.
Set $i = ( ( 2d )! , p-1 )$ and $k_i$ the subfield of $k ( \zeta_p )$ of degree $i$ over $k$.
If for every $P \in \A[p]$ of order $p$ the field $k ( P ) \cap k ( \zeta_p )$ strictly contains $k_i$, then $\Sh ( \A / k )$ is $p$-divisible in $H^1 ( k , \A )$,
\end{cor}  
     
\begin{cor}\label{cor7}
Let $ \A $ be a polarized abelian surface defined over $ k $.
For every prime number $p > 3840$ such that $k \cap \Q ( \zeta_p ) = k$ and not dividing the degree of the polarization, then if $\Sh ( \A / k )$ is not $p$-divisible over $H^1 ( k , \A )$, there exists a finite extension $\widetilde{k}$ of $ k $ of degree $\leq 6$ such that $ \A^t $ is $\widetilde{k}$-isogenous to an abelian surface with a torsion point of order $ p $ defined over $ \widetilde{k} $.
\end{cor}

\end{document}